\input amstex 
\documentstyle{amsppt}
\input bull-ppt
\keyedby{bull390e/mhm}


\topmatter
\cvol{28}
\cvolyear{1993}
\cmonth{April}
\cyear{1993}
\cvolno{2}
\cpgs{288-305}
\title Wavelet transforms versus Fourier transforms 
\endtitle
\author Gilbert Strang \endauthor
\address Department of Mathematics,
Massachusetts Institute of Technology,
Cambridge, Massachusetts 02139\endaddress
\ml gs\@math.mit.edu\endml
\keywords Wavelets, Fourier transform, dilation,
orthogonal basis\endkeywords  
\subjclass Primary 42A06, 41A05, 65D05\endsubjclass
\thanks I am grateful to the National Science Foundation 
(DMS 90-06220) 
for their support\endthanks
\abstract This note is a very basic introduction to 
wavelets.  It
starts with an orthogonal basis of piecewise constant 
functions,
constructed by dilation and translation.  The ``wavelet 
transform''
maps each $f(x)$ to its coefficients with respect to this 
basis.  The
mathematics is simple and the transform is fast (faster 
than the Fast
Fourier Transform, which we briefly explain), 
but approximation by piecewise constants is poor.  To 
improve this
first wavelet, we are led to dilation equations and their 
unusual
solutions.  Higher-order wavelets are constructed, and it is
surprisingly quick to compute with them --- always 
indirectly and
recursively.

We comment informally on the contest between these 
transforms in
signal processing, especially for video and image 
compression
(including high-definition television).
So far the Fourier Transform --- or its 8 by 8 windowed 
version, the Discrete Cosine Transform --- is often chosen. 
But wavelets are already competitive, and they are ahead 
for fingerprints.
We present a sample of this developing theory.\endabstract 
\date March 20, 1992 and, in revised form, November 30, 
1992\enddate
\endtopmatter

\document

\heading 1. The Haar wavelet \endheading

To explain wavelets we start with an example. It has
every property we hope for, except one.  If that one 
defect is
accepted, the construction is simple and the computations 
are fast.
By trying to remove the defect, we are led to dilation 
equations and
recursively defined functions and a small world of 
fascinating new
problems --- many still unsolved.  A sensible person would 
stop after
the first wavelet, but fortunately mathematics goes on.

The basic example is easier to draw than to describe:

\midspace{6.5pc}
\caption{{\smc Figure} 1. 
{\rm Scaling function $\phi (x)$, wavelet $W(x)$,
and the next level of detail.}}

Already you see the two essential operations:  {\it 
translation}
and {\it dilation}.  The step from $W(2x)$ to $W(2x-1)$ is
translation.  The step from $W(x)$ to $W(2x)$ is dilation. 
 Starting
from a single function, the graphs are shifted and 
compressed.
The next level contains\ $W(4x)$, $W(4x-1)$, $W(4x-2)$, 
$W(4x-3)$. Each
is supported on an interval of length $\frac{1}{4}$.  In 
the end we have
Haar's infinite family of functions:
$$W_{jk}(x) = W(2^{j}x-k) \quad (\text{together with}\
\phi(x)).$$
When the range of indices is $j \geq 0$ and $0 \leq k < 
2^{j}$, these
functions form a remarkable basis for $L^{2}[0,1]$.  We 
extend it
below to a basis for $L^{2}$({\bf R}).  

The four functions in Figure~1 are piecewise constant.  
Every function that is
constant on each quarter-interval is a combination of 
these four.  Moreover,
the inner product  $\int \phi (x) \, W(x) \, dx$ is zero 
--- and so are the
other inner products.  This property extends to all $j$ 
and $k$: {\it The
translations and dilations} {\it of}  $W$ {\it are 
mutually orthogonal.}  We
accept  this as the definition of a wavelet, although 
variations are definitely
useful in practice.  The goal looks easy enough, but the 
example is deceptively
simple.  

This orthogonal Haar basis is not a recent invention [1].  
It is
reminiscent of the Walsh basis in [2] --- but the 
difference is
important.\footnote"$^\dagger$"{Rademacher was first to 
propose an orthogonal
family of $\pm 1$ functions; it was not complete.  After 
Walsh constructed a
complete set, Rademacher's Part II was regrettably 
unpublished and seems to be
lost (but Schur saw it).}  For Walsh and Hadamard, the 
last two basis
functions are changed to $W(2x) \pm W(2x-1)$.  All of their
``binary sinusoids'' are
supported on the 
whole interval $0 \leq x \leq 1$.  This global support is 
the one drawback to
sines and 
cosines; otherwise, Fourier is virtually unbeatable.
To represent a local function, vanishing outside a short
interval of space or time, a global basis requires extreme
cancellation.  Reasonable accuracy needs many terms of the 
Fourier series.
{\it Wavelets give a local basis}.

You see the consequences.  If the signal $f(x)$ disappears 
after $x =
\frac{1}{4}$, only a quarter of the later basis functions 
are 
involved.  The wavelet 
expansion directly reflects the properties of $f$ in 
physical space,
while the Fourier expansion is perfect in frequency space. 
 Earlier
attempts at a ``windowed Fourier transform'' were ad hoc 
--- wavelets
are a systematic construction of a local basis.

The great value of orthogonality is to make expansion
coefficients easy to compute.  Suppose the values of 
$f(x)$, constant
on four quarter-intervals, are $9,1,2,0$.  Its Haar 
wavelet expansion
expresses this vector $y$ as a combination of the basis 
functions:
$$\bmatrix 9 \\ 1 \\ 2 \\ 0 \endbmatrix =
3 \bmatrix 1 \\ 1 \\ 1 \\ 1\endbmatrix +
2 \bmatrix \format\r\\ \ 1 \\ \ 1 \\ -1 \\ -1\endbmatrix +
4\bmatrix \format\r\\ \ 1 \\ -1 \\ \ 0 \\ \ 0\endbmatrix +
\bmatrix \format\r\\ \ 0 \\ \ 0 \\ \ 1 \\ -1\endbmatrix.
$$
The wavelet coefficients $b_{jk}$ are $3,2,4,1$; they form 
the
wavelet transform of $f$.  The connection between the 
vectors $y$
and $b$ is the matrix $W_{4}$, in
whose orthogonal columns you recognize the graphs of 
Figure~1:
$$y = W_{4}b\ \ \ \text{is}\ \ \
\bmatrix 9 \\ 1\\ 2\\ 0\endbmatrix =
\bmatrix \format\r&\quad\r&\quad\r&\quad\r \\
1 & 1 & 1 & 0 \\
\ 1 & \ 1 & -1 & \ 0 \\
\ 1 & -1 & \ 0 & \ 1 \\
\ 1 & -1 & \ 0 & -1 \endbmatrix
\bmatrix 3 \\ 2 \\ 4 \\ 1\endbmatrix .$$
This is exactly comparable to the {\it Discrete Fourier
Transform}, in which  $f(x) = \sum a_{k} \, e^{ikx}$
stops after four terms.  Now
the vector $y$ contains the values of $f$ at four points:
$$y = F_{4}a\ \ \ \text{is}\ \ \ \bmatrix f(0\pi/2) \\ 
f(1\pi/2) \\ f(2\pi
/2) \\ f(3\pi /2) \endbmatrix = \bmatrix 1& 1 & 1 & 1 \\ 1 
& i
&i^{2} & i^{3} \\ 1 &i^{2} & i^{4} & i^{6} \\ 1 &i^{3} & 
i^{6}
& i^{9}\endbmatrix \bmatrix a_{0} \\ a_{1} \\ a_{2} \\
a_{3}\endbmatrix .$$
This Fourier matrix also has orthogonal columns.  The $n$ 
by $n$
matrix $F_{n}$ follows the same pattern, with $\omega = 
e^{2\pi
i/n}$ in place of $i = e^{2\pi i/4}$.  Multiplied by
$1/\sqrt{n}$ to give orthonormal columns, it is the
most important of all unitary matrices.  The wavelet matrix
sometimes offers modest competition.

To invert a real orthogonal matrix we transpose it.  To 
invert a
unitary matrix, transpose its complex conjugate.  After 
accounting
for the factors that enter when columns are not unit 
vectors, the
inverse matrices are
$$W^{-1}_{4} = \frac{1}{4} \bmatrix
\format\r&\quad\r&\quad\r&\quad\r \\
1 & 1 & 1 & 1 \\
\ 1 & \ 1 & -1 & -1 \\
\ 2 & -2 & \ 0 & \ 0 \\
\ 0 & \ 0 & \ 2 & -2 \endbmatrix
\ \ \ \ \text{and} \ \ \
F^{-1}_{4} = \frac{1}{4} \bmatrix 1 & 1 & 1 & 1 \\
1 & (-i) & (-i)^{2} & (-i)^{3} \\ 1 & (-i)^{2} & (-i)^{4} 
& (-i)^{6}
\\ 1 & (-i)^{3} &  (-i)^{6} & (-i)^{9} \endbmatrix.$$
The essential point is that the inverse matrices have the 
same form
as the originals.  If we can transform quickly, we can 
invert quickly ---
between coefficients and function values.  The Fourier 
coefficients come from
values at $n$ points.  The Haar coefficients come from 
values on $n$
subintervals.

\heading 2. Fast Fourier Transform and Fast Wavelet 
Transform \endheading

The Fourier matrix is full --- it has no zero entries. 
Multiplication of
$F_{n}$ times a vector $a$, done directly, requires 
$n^{2}$ separate
multiplications.  We are evaluating an $n$-term Fourier 
series at $n$ points. 
The series is $\sum^{n-1}_{0} a_{k} \, e^{ikx}$, and the 
points are $x = 2\pi
j/n$.

The wavelet matrix is sparse --- many of its entries are 
zero.
Taken together, the third
and fourth columns of $W$ fill a single column; the fifth,
sixth, seventh, and eighth columns would fill one more 
column.
With $n = 2^{\ell}$, we fill only $\ell + 1$ columns.  The 
total
number of nonzero entries in
$W_{n}$ is $n(\ell + 1)$.  This already shows the effect 
of a more
local basis.  Multiplication of $W_{n}$ times a vector 
$b$, done 
directly, requires only $n(\log_{2} n + 1 )$ separate
multiplications.

Both of these matrix multiplications can be made faster.  
For
$F_{n}a$, this is achieved by the Fast Fourier Transform 
--- the most
valuable numerical algorithm in our lifetime.  It changes 
$n^{2}$ to
$\frac{1}{2} n \log_{2}n$ by a reorganization of the steps 
--- which
is simply a factorization of the Fourier matrix.  A typical
calculation with $n = 2^{10}$ changes $(1024)(1024)$
multiplications to $(5)(1024)$.  This saving by a factor 
greater than
$200$ is genuine.  The result is that the FFT has 
revolutionized
signal processing.  Whole industries are changed from slow 
to fast by
this one idea --- which is pure mathematics.

The wavelet matrix $W_{n}$ also allows a neat 
factorization into very
sparse matrices.  The operation count drops from $O(n \log 
n)$ all
the way to $O(n)$.  For our piecewise constant wavelet the 
only
operations are add and subtract; in fact, $W_{2}$ is the 
same as
$F_{2}$.  Both fast transforms have  $\ell = \log_{2}n$ 
steps, in the
passage from $n$ down to 1.  For the FFT, each step  
requires $\frac{1}{2}n$
multiplications (as shown below).  For the Fast Wavelet 
Transform,
the cost of each successive step is cut in half.  It is a 
beautiful ``pyramid
scheme'' created by Burt and Adelson and Mallat and 
others. The total cost has
a factor $1 + \frac{1}{2} + \frac{1}{4} 
+ \cdots$ that stays below $2$.  This is why the final 
outcome for
the FWT is $O(n)$ without the logarithm $\ell$.

The matrix factorizations are so simple, especially for $n 
= 4$, that it
seems worthwhile to display them.  The FFT has two copies 
of the
half-size transform $F_{2}$ in the middle:
$$F_{4} = \bmatrix \format \r &\quad\r &\quad \r &\quad \r 
 \\1 & \ &  1 & \
\\  \ & 1 &\ & i \\ 1 &\ &-1 &\ \\ \
& 1 &\ & -i \endbmatrix  \bmatrix 1 & 1 &\ & \ \\ 1 & 
i^{2} &\  & \
\\ \ & \  & 1 & 1 \\ \ & \ & 1 & i^{2}\endbmatrix \bmatrix 
1 &\ &\ &\
\\ \ & \ & 1 &\ \\ \ & 1 &\ &\  \\ \ & \ & \ &1\endbmatrix .
\tag{1}$$
The permutation on the right puts the even $a$\<'s
($a_{0}$ and $a_{2}$) ahead of the odd $a$\<'s ($a_{1}$ 
and $a_{3}$).
Then come {\it separate half-size transforms} on the evens 
and odds.
The matrix at the left combines these two half-size 
outputs in a
way that produces the correct full-size answer.  By 
multiplying
those three matrices we recover $F_{4}$.

The factorization of $W_{4}$ is a little different:
$$W_{4} = \bmatrix 1 & \ 1 &\ &\ \\ 1 & -1 &\ &\ \\ \ & \ 
& 1 & \ 1 \\ \
& \ & 1 & -1\endbmatrix  \bmatrix 1 & \ &\ &\ \\ \ & \ & 1 
&\ \\ \ &
1 &\ &\ \\ \ & \ &\ & 1\endbmatrix \bmatrix  \format \r 
&\quad \r &\quad \r
&\quad  \l \\ 1 & 1 &\ &\ \\ 1 &-1 &\
&\ \\ \ &\ & 1 &\ \\ \ &\ &\ &1\endbmatrix . \tag{2}$$
At the next level of detail (for $W_{8}$), the same $2$ by 
$2$ matrix
appears four times in the left factor.  The permutation 
matrix puts
columns $0,2,4,6$ of that factor ahead of $1,3,5,7$.  The 
third factor has
$W_{4}$ in one corner and $I_{4}$ in the other
corner  (just as $W_{4}$ above ends with $W_{2}$ and 
$I_{2}$ --- this
factorization is the matrix form of the pyramid 
algorithm).  It
is the identity matrices $I_{4}$ and $I_{2}$  that save
multiplications.  Altogether $W_{2}$ 
appears $4$ times at the left of $W_{8}$, then $2$ times 
at the left of
$W_{4}$, and then once at the right.  The multiplication 
count from these $n-1$
small matrices is $O(n)$ --- the
Holy Grail of complexity theory.

Walsh would have another copy of the $2$ by $2$ matrix in 
the last
corner, instead of $I_{2}$.  Now the product
has orthogonal columns with all entries $\pm 1$ --- the 
Walsh basis.  
Allowing $W_{2}$ or $I_{2}$, $W_{4}$ or $I_{4}$, $W_{8}$ 
or $I_{8}$,
$\dots$ in the third factors, the matrix products exhibit 
a whole family of
orthogonal bases.  This is a  {\it wavelet
packet}, with great flexibility.  Then a ``best basis'' 
algorithm aims for a
choice 
that concentrates most of $f$ into a few basis vectors. 
That is the goal ---
to compress information.

The same principle of factorization applies for any power 
of $2$, say
$n = 1024$.  For
Fourier, the entries of $F$ are
powers of $\omega = e^{2\pi i/1024}$.  The row and column 
indices
go from $0$ to $1023$ instead of $1$ to $1024$.  The 
zeroth row and
column are filled with  $\omega^{0} = 1$.  {\it The
entry in row} $j$, {\it column} $k$ {\it of} $F$
{\it is} $\omega^{jk}$. This is the term $e^{ikx}$ 
evaluated at $x = 2\pi
j/1024$.  The multiplication $F_{1024}a$  computes the 
series $\sum
a_{k} \, \omega^{jk}$ for $j
= 0$ to $1023$. 

The key to the matrix factorization is just this.  
Squaring the
$1024\text{th}$ root
of unity gives the $512$\<th root: $(\omega^{2})^{512} = 
1$.  This was the
reason behind 
the middle factor in (1), where $i$ is the fourth root and 
$i^{2}$ is the
square root.  It is the essential link
between $F_{1024}$ and $F_{512}$.  The first stage of the 
FFT is the
great factorization rediscovered by Cooley and Tukey (and 
described in
1805 by Gauss):
$$F_{1024} = \bmatrix I_{512} & \ D_{512} \\ I_{512} &
-D_{512}\endbmatrix  \bmatrix F_{512} &\ \\ \ & 
F_{512}\endbmatrix
\bmatrix \text{even-odd} \\ \text{shuffle}\endbmatrix .
\tag{3}$$
$I_{512}$ is the identity matrix.  $D_{512}$ is the 
diagonal matrix
with entries $(1,\omega, 
\dots ,\omega^{511})$, requiring about $512$
multiplications.  The two copies of $F_{512}$ in the 
middle give a
matrix only {\it half full} compared to $F_{1024}$ --- 
here is
the crucial saving.  The shuffle separates the incoming
vector $a$ into $(a_{0},a_{2},\dots ,a_{1022})$ with even
indices and the odd part $(a_{1},a_{3},\dots ,a_{1023})$.

Equation (3) is an imitation of equation (1), eight levels 
higher.  Both are
easily verified.  Computational linear algebra has become 
a world of matrix
factorizations, and this one is outstanding.

You have anticipated what comes next.  Each $F_{512}$ is 
reduced
in the same way to two half-size transforms $F = F_{256}$. 
 The work is
cut in half 
again, except for an additional $512$ multiplications from 
the
diagonal matrices $D = D_{256}$: 
$$\bmatrix 
{}\\
F_{512}  \\
& F_{512} \\
{}
\endbmatrix =
\bmatrix 
\format \r &\ \ \r & \r &\ \ \r \\
I & D & \ & \ \\
I & -D & \ & \\
\ & \ & I & D \\
\ & \ & I & -D \endbmatrix
\bmatrix F  \\
& F  \\
& & F  \\
& & & F \endbmatrix
\bmatrix \text{even-odd gives} \\
0\ \text{and}\ 2\mod 4 \\
\text{even-odd gives} \\
1\ \text{and}\ 3\mod 4 \endbmatrix\.\tag 4$$  
For $n = 1024$  there are $\ell = 10$ levels, and each level
has $\frac{1}{2}n = 512$ multiplications  from the first 
factor --- to
reassemble  the half-size outputs
from the level below.  Those $D$\<'s yield
the final count $\frac{1}{2}n\ell$.  

In practice, $\ell = \log_{2}n$ is controlled by splitting
the signal into smaller blocks. With $n = 8$, the scale 
length of the transform
is closer to the scale length of most images. 
This is the {\it short time} Fourier transform, which is the
transform of a ``windowed'' function $wf$.  The multiplier 
$w$
is the characteristic function of the window.  (Smoothing is
necessary!  Otherwise this
blocking of the image can be visually unacceptable.  The 
ridges of
fingerprints are broken up very badly, and windowing was
unsuccessful in tests by the FBI.) 
In other applications the implementation may favor the FFT 
---
theoretical complexity is rarely the whole story.

A more gradual exposition of the Fourier matrix and the 
FFT is in
the monographs [3, 4] and the textbooks [5, 6] --- and in 
many other
sources [see 7].  (In the lower level text [8], it is 
intended 
more for reference
 than for teaching.  On the other hand, this is just a 
matrix--vector
multiplication!)  FFT codes are
freely available on {\it netlib}, and generally each 
machine has its
own special software.

For higher-order wavelets, the FWT still involves many 
copies of a
single small matrix.  The entries of this matrix are 
coefficients
$c_{k}$ from the ``dilation equation''. We
move from fast
algorithms to a quite different part of mathematics --- 
with the goal
of constructing new orthogonal bases.  The basis functions 
are 
unusual, for a good reason.

\heading 3. Wavelets by multiresolution analysis \endheading

The defect in piecewise constant wavelets is that they are
very poor at approximation.  Representing a smooth
function requires many pieces.  For wavelets
this means many levels --- the number $2^{j}$ must be 
large for an
acceptable accuracy.  It is similar to the rectangle rule 
for integration,
or Euler's method for a differential equation, or forward 
differences
$\Delta y /\Delta x$ as estimates of $dy/dx$. Each is a 
simple and
natural first approach, but inadequate in the end.  
Through all of
scientific computing 
runs this common theme:  Increase the accuracy at least to
second order.  What this means is: {\it Get the
linear term right}.  

For integration, we move to the trapezoidal rule and 
midpoint rule.  For derivatives, second-order accuracy 
comes with centered
differences.  The whole point of Newton's
method for solving $f(x) = 0$ is to follow the tangent 
line.  All
these are exact when $f$ is linear.  For wavelets to be 
accurate, 
$W(x)$ and $\phi (x)$ need the same improvement.  Every 
$ax+b$ must be a
linear combination of translates.

Piecewise polynomials
(splines and finite elements) are often based on the 
``hat'' function --- the
integral of Haar's $W(x)$.  But this
 piecewise linear function does not produce
orthogonal wavelets with a local basis.  The requirement of
{\it orthogonality to dilations} conflicts strongly with 
the demand for
compact support --- so much so that it was originally 
doubted whether one
function could satisfy both requirements and still produce 
$ax+b$.  It was
the achievement of
Ingrid Daubechies [9] to construct such a function.
 
We now outline the construction of wavelets.  The reader 
will understand
that we only touch  on parts of the theory 
and on selected applications.  An excellent account of the 
history is in [10].
Meyer and Lemari\'{e} describe the earliest wavelets 
(including Gabor's). Then
comes the beautiful pattern of multiresolution analysis 
uncovered by Mallat --- which is hidden by the simplicity of
the Haar basis.  Mallat's analysis found expression
in the Daubechies wavelets.

Begin on the interval $[0,1]$. The space $V_0$ spanned by 
$\phi(x)$ is 
orthogonal to the space $W_0$ spanned by $W(x)$. Their sum 
$V_1=V_0\oplus W_0$
consists of all piecewise constant functions on 
half-intervals.  A different
basis for $V_1$ is $\phi(2x)=\frac 12 (\phi(x)+W(x))$ and 
$\phi(2x-1)=\frac 
12 (\phi(x)-W(x))$. Notice especially that $ V_0\subset  
V_1$. 
{\it The function
$\phi(x)$ is a combination of $\phi(2x)$ and $\phi(2x-1)$}. 
This is the dilation 
equation, for Haar's example.

Now extend that pattern to the spaces $V_j$ and $W_j$ of 
dimension $2^j$:
$$\align
V_{j} &= \text{span of the translates} \ \phi (2^{j}x - 
k)\ \text{for fixed}\
j ,\\
W_{j} &= \text{span of the wavelets}\ W(2^{j}x - k)\ 
\text{for fixed}\ j.
\endalign
$$
The next space $V_{2}$ is spanned by $\phi (4x)$, $\phi 
(4x-1)$,
$\phi(4x - 2)$, $\phi (4x-3)$.  It contains all piecewise 
constant
functions on quarter-intervals.  That space was also 
spanned by the four
functions $\phi (x)$, $W(x)$, $W(2x)$, $W(2x-1)$ at the 
start of this paper.
Therefore, $V_{2}$ decomposes into $V_{1}$ and $W_{1}$ 
just as $V_{1}$
decomposes into $V_{0}$ and $W_{0}$:
$$V_{2} = V_{1} \oplus W_{1} = V_{0} \oplus W_{0} \oplus 
W_{1}. \tag{5}$$
At every level, {\it the wavelet space} $W_{j}$ {\it is the
``difference''} 
{\it between} $V_{j+1}$ {\it and} $V_{j}$:
$$V_{j+1} = V_{j} \oplus W_{j} = V_{0} \oplus W_{0} \oplus 
\cdots \oplus
W_{j}. \tag{6}$$
The translates of wavelets on the right are also 
translates of scaling
functions on the left.  
For the construction of
wavelets, this offers a totally different approach.  
Instead of creating
$W(x)$ and the spaces $W_{j}$, we can create $\phi (x)$ 
and the spaces
$V_{j}$.  It is a choice between the terms $W_{j}$ of an 
infinite series or
their 
partial sums $V_{j}$.  Historically the constructions 
began with $W(x)$. Today
the 
constructions begin with $\phi (x)$.  It has proved easier 
to work with sums
than differences.

A first step is to change from $[0,1]$ to the whole line 
{\bf R}.  The
translation index $k$ is unrestricted. The subspaces
$V_{j}$ and $W_{j}$ are
infinite-dimensional ($L^{2}$ closures of translates).  
One basis for
$L^2(\bold R)$ consists of $\phi(x-k)$ and 
$W_{j\,k}(x)=W(2^jx-k)$ with $j\ge 0,
k\in \bold Z$. Another basis contains all $W_{j\,k}$ with 
$j,k\in \bold Z$. Then
the 
dilation
index $j$ is also unrestricted --- for $j = -1$
the functions $\phi (2^{-1}x-k)$ are
constant on intervals of length $2$.  
The decomposition into $V_{j} \oplus W_{j}$
continues to hold!  The sequence of closed subspaces $V_{j}$
has the following basic properties for $-\infty < j < 
\infty$:
$$\align
&V_{j} \subset V_{j+1}\ \ \text{and}\ \ \bigcap V_{j} = \{ 
0\}\ \ \text{and}\ \
\bigcup V_{j}\ \text{is dense in} \ L^{2}({\bold R}) ;\\
&f(x) \ \text{is in}\ V_{j}\ \text{if and only if}\ f(2x)\ 
\text{is in}\
V_{j+1} ;\\
&V_{0}\ \text{has an orthogonal basis of translates}\ \phi 
(x-k),\ k \in
{\bold Z}.
\endalign
$$
These properties yield a {\rm ``}{\it multiresolution 
analysis}\/{\rm ''} 
--- the
pattern that other wavelets will follow.  $V_{j}$ will be 
spanned by $\phi
(2^{j}x - k)$.  $W_{j}$ will be its orthogonal complement 
in $V_{j+1}$.
Mallat proved, under mild hypotheses, that $W_{j}$ is also 
spanned
by translates [11]; these are the wavelets.

Dilation is built into multiresolution analysis by the 
property that
$f(x) \in V_{j} \Leftrightarrow f(2x) \in V_{j+1}$.
This applies in particular to $\phi (x)$.
It must be a combination of translates of $\phi (2x)$. 
That is the hidden
pattern, which has become central to this subject. We have 
reached
the dilation equation.

\heading 4. The dilation equation \endheading

In the words of [10], ``{\it la perspective est 
compl\`{e}tement
chang\'{e}e}.'' The 
construction of wavelets now begins with the scaling 
function $\phi$.  
The dilation
equation (or refinement equation or two-scale difference 
equation)
connects $\phi (x)$ to translates of $\phi (2x)$:
$$\phi (x) = \sum^{N}_{k=0} c_{k}\ \phi (2x-k). \tag{7}$$
The coefficients for Haar are $c_{0} = c_{1} = 1$.  The 
box function $\phi$ is
the sum of two half-width boxes. That is equation (7).  
Then $W$ is a
combination of the same translates (because $W_{0} \subset 
V_{1}$).  The
coefficients for $W = \phi (2x) - \phi (2x-1)$ are $1$ and 
$-1$. It is
absolutely remarkable that $W$ uses the same
coefficients as $\phi$, but in reverse order and with 
alternating signs:
$$W(x) = \sum^{1}_{1-N} \ (-1)^{k}\ c_{1-k} \ \phi(2x-k). 
\tag{8}$$
This construction makes $W$ orthogonal to $\phi$ and its 
translates.
(For those translates to be orthogonal to each other, see
below.)  The key is that every vector 
$c_{0},c_{1},c_{2},c_{3}$ is
automatically orthogonal to $c_{3},-c_{2},c_{1},-c_{0}$ 
and all even
translates like $0,0,c_{3},-c_{2}$.

When $N$ is odd, $c_{1-k}$ can be replaced in (8) by 
$c_{N-k}$.  This
shift by $N-1$ is even.  Then the sum goes from $0$ to $N$ 
and $W(x)$
looks especially attractive.

{\it Everything hinges on the $c$\<{\rm '}s\/}. They 
dominate all that follows.
They 
determine (and are determined by) $\phi$, they determine 
$W$, and they go into
the matrix factorization (2).  In the applications, 
convolution with $\phi$ is
an averaging operator --- it produces smooth 
functions (and a blurred picture). Convolution with $W$ is 
a differencing 
operator, which picks out
details. 

The convolution of the box with itself is the piecewise 
linear hat function
---  equal to $1$ at $x = 1$ and supported on the interval 
$[0,2]$.  It
satisfies the dilation equation with
$c_{0} = \frac{1}{2}$, $c_{1} = 1$, $c_{2} = \frac{1}{2}$.
But there is a condition on the $c$\<'s in order that
the wavelet basis $W(2^{j}x-k)$ shall be orthogonal.  The 
three
coefficients $\frac{1}{2},1,\frac{1}{2}$ do not satisfy 
that condition.
Daubechies found the unique
$c_{0},c_{1},c_{2},c_{3}$ (four coefficients are 
necessary) to give
orthogonality plus second-order
approximation. Then the question becomes: {\it How to 
solve} {\it the
dilation equation}? 

\medskip

{\it Note added in proof.} A new construction has just 
appeared
that uses {\it two} scaling functions $\phi_{i}$ and 
wavelets
$W_{i}$.  Their translates are still orthogonal [38].  The 
combination
$\phi_{1}(x) + \phi_{1}(x-1) + \phi_{2}(x)$ is the hat 
function, so
second-order accuracy is achieved.  The remarkable 
property is that
these are ``{\it short functions}'': $\phi_{1}$ is 
supported on
$[0,1]$ and $\phi_{2}$ on $[0,2]$.  They satisfy a matrix 
dilation
equation.

These short wavelets open new possibilities for 
application, since
the greatest difficulties are always at boundaries.  The 
success of
the finite element method is largely based on the very 
local character
of its basis functions.  Splines have longer support (and 
more smoothness),
wavelets have even longer support (and orthogonality).  The
translates of a long basis function overrun the boundary.   

\medskip

There are two principal methods to solve dilation 
equations.  One is by Fourier transform, the
other is by matrix products.  Both give $\phi$ as a limit, 
not as an
explicit function. We never
discover the exact value $\phi(\sqrt{2})$.  It is amazing to
compute with a function we do not know --- but the 
applications only require
the 
$c$\<'s.  When complicated functions come from a simple 
rule, we know
from increasing experience what to do: Stay with the 
simple rule.

\subheading{Solution of the dilation equation by Fourier 
transform} Without
the ``2'' we would have an 
ordinary difference equation --- entirely familiar.  The 
presence of two
scales, $x$ and $2x$, is the problem.  A warning comes 
from Weierstrass
and de Rham and Takagi --- their nowhere differentiable 
functions are
all built on multiple scales like $\sum a^{n} \cos 
(b^{n}x)$.  The
Fourier transform easily handles translation by $k$ in 
equation (7),
but $2x$ in physical
space becomes $\xi /2$ in frequency space:
$${\hat \phi}(\xi) = \frac{1}{2} \sum c_{k} \, e^{ik\xi 
/2} \,
{\hat \phi}\left(\frac{\xi}{2}\right) =
P\left(\frac{\xi}{2}\right) \, {\hat 
\phi}\left(\frac{\xi}{2}\right). \tag{9}$$
The ``symbol'' is
$P(\xi) = \frac{1}{2} \sum c_{k} \, e^{ik\xi}$.  With $\xi =
0$ in (9) we find $P(0) = 1$ or $\sum c_{k} = 2$ --- the 
first requirement
on the $c$\<'s.  This allows us to look for a solution 
normalized by ${\hat
\phi}(0) = \int \phi (x) \, dx = 1$.  It does not ensure 
that we find a
$\phi$ that is continuous or even in $L^{1}$.
 What we do find is an infinite
product, by recursion from $\xi /2$ to $\xi /4$ and onward:
$${\hat \phi}(\xi) = P\left(\frac{\xi}{2}\right)\, {\hat 
\phi}
\left(\frac{\xi}{2}\right) =
P\left(\frac{\xi}{2}\right) \, P\left(\frac{\xi}{4}\right) 
\, {\hat \phi}\left(\frac{\xi}{4}\right) =
\cdots =
\prod^{\infty}_{j=1} P\left(\frac{\xi}{2^{j}}\right).$$
This solution $\phi$ may be only a distribution.
Its smoothness becomes clearer by matrix
methods. 

\subheading{Solution by matrix products \cite {12, 13}} 
When $\phi$ is
known at the integers, the dilation equation gives $\phi$ 
at 
half-integers such as $x = \frac{3}{2}$.  Since $2x-k$ is 
an integer, we just
evaluate  $\sum c_{k} \phi (2x-k)$.  Then the equation 
gives $\phi$ at quarter-integers as combinations of $\phi$ 
at half-integers.
 The combinations are built into the entries of two 
matrices $A$ and $B$, and
the recursion
is taking their products.

To start we need $\phi$ at the integers.  With $N = 3$, for
example, set $x = 1$ and $x = 2$ in the dilation equation:
$$\aligned
\phi(1) &= c_{1}\,\phi(1) + c_{0} \,\phi(2) ,\\
\phi(2) &= c_{3}\,\phi(1) + c_{2} \,\phi(2). 
\endaligned \tag{10}$$
Impose the conditions $c_{1} + c_{3} = 1$ and $c_{0} + 
c_{2} = 1$.  Then the
$2$~by~$2$ matrix in (10), formed from these $c$\<'s, has 
$\lambda = 1$
as an eigenvalue.  The eigenvector
is $(\phi (1) , \phi (2))$.  It follows from (7) that 
$\phi$ will
vanish outside $0 \leq x \leq N$.  

To see the step from integers to half-integers in matrix 
form,
convert the scalar dilation equation to a
first-order equation for the vector $v(x)$:
$$v(x) = \bmatrix \phi(x) \\ \phi(x+1) \\ \phi(x+
2)\endbmatrix, \qquad 
A =
\bmatrix c_{0} & 0 & 0 \\ c_{2} & c_{1} & c_{0} \\ 0 &  
c_{3} &
c_{2}\endbmatrix, \qquad 
B = \bmatrix c_{1} & c_{0} & 0 \\ c_{3} &
c_{2} & c_{1} 
\\ 0 & 0 & c_{3} \endbmatrix.$$
The equation turns out to be $v(x) = Av(2x)$ for $0 \leq x 
\leq
\frac{1}{2}$ and 
$v(x) = Bv(2x-1)$ for
$\frac{1}{2} \leq x \leq 1$.  By recursion this yields $v$
at any dyadic point --- whose binary expansion is finite.  
Each $0$
or $1$ in the expansion decides between $A$ and $B$.  For 
example
$$v(.01001) = (ABAAB)v(0). \tag{11}$$

Important: {\it The matrix} $B$ {\it has entries} 
$c_{2i-j}$.  So
does $A$, when the indexing starts with $i = j = 0$.  The
dilation equation itself is $\phi = C\phi$, with an 
operator $C$ of this new 
kind.  Without the $2$ it would be a Toeplitz operator, 
constant along 
each diagonal, but now every
other row is removed.  Engineers call it ``convolution 
followed by
decimation''.  (The word {\it downsampling} is also used 
--- possibly a
euphemism for decimation.)  Note that the derivative of 
the dilation
equation is $\phi^{\prime} = 
2C\phi^{\prime}$.  Successive derivatives introduce powers 
of $2$.
The eigenvalues of these operators $C$ are $1,\frac{1}{2},
\frac{1}{4},\dots$, until $\phi^{(n)}$ is not defined in 
the space
at hand.  The sum condition $\sum c_{\text{even}} = \sum
c_{\text{odd}} = 1$ is always imposed --- it assures in 
Condition
$\text{A}_{1}$ below that we have 
first-order approximation at least.

When $x$ is not a dyadic point $p/2^{n}$, the recursion in 
(11) does not
terminate.  The
binary expansion $x = .0100101\dots$ corresponds to an 
infinite product
$ABAABAB\dots$.  The convergence of such a product is
by no means assured.  It is a major problem to
find a direct test on the $c$\<'s that is equivalent to 
convergence --- for
matrix products {\it in every order}.  We briefly describe 
what
is known for arbitrary $A$ and $B$.

For a single matrix $A$, the growth of the powers $A^{n}$ 
is governed
by the spectral radius $\rho(A) = \max |\lambda_{i}|$.  
Any norm of
$A^{n}$ is roughly the $n$\<th power of this largest 
eigenvalue.
Taking $n$\<th roots makes this precise:
$$\lim_{n \rightarrow \infty}\  \Vert A^{n}\Vert^{1/n} = 
\rho(A).$$
The powers approach zero if and only if $\rho (A) < 1$.

For two or more matrices, the same process produces the 
{\it joint
spectral radius} [14].  The powers $A^{n}$ are replaced by 
products $\Pi_{n}$
of $n\ A$\<'s and $B$\<'s.  The maximum of $\| \Pi_{n}\|$, 
allowing products in
all orders, is still submultiplicative. The limit of 
$n$\<th roots (also the
infimum) is the joint spectral radius:
$$\lim_{n \rightarrow \infty}\ (\max \Vert 
\Pi_{n}\Vert)^{1/n} =
\rho(A,B). \tag{12}$$
The difficulty is not to define $\rho(A,B)$ but to compute 
it.  For
symmetric or normal or commuting or upper triangular  
matrices it is
the larger of $\rho(A)$ and
$\rho(B)$.  Otherwise eigenvalues of
products are not controlled by products of eigenvalues.  
An example with
zero eigenvalues, $\rho (A) = 0 = \rho (B)$, is
$$A = \left[ \matrix 0 & 2 \\ 0 & 0\endmatrix \right] 
,\qquad B =
\left[ \matrix 0 & 0
\\ 2 & 0\endmatrix \right] ,\qquad
AB = \left[ \matrix 4 & 0 \\ 0 & 0
\endmatrix \right].$$
In this case  $\rho
(A,B) = \| AB \|^{1/2} = 2$.  The product $ABABAB\dots$ 
diverges.  
In general
$\rho$ is a function of the matrix entries, bounded above 
by norms and
below by eigenvalues.  Since one possible infinite product 
is a
repetition of any particular $\Pi_{n}$ (in the example it 
was $AB$), the
spectral radius of that
single matrix gives a lower bound on the joint radius: 
$$(\rho(\Pi_{n}))^{1/n} \leq \rho(A,B). $$
A beautiful theorem of Berger and Wang [15] asserts that 
these
eigenvalues of products yield the same limit  (now a 
supremum)
that was approached by norms:
$$\underset{n \rightarrow \infty}\to{\lim \sup}\ ( \max 
\rho(\Pi_{n}))^{1/n} =
\rho(A,B). \tag{13}$$
It is conjectured by Lagarias and Wang that equality is 
reached at a
finite product 
$\Pi_{n}$.  Heil and the author noticed a corollary of the
Berger-Wang theorem: $\rho$ is a continuous function of 
$A$ and $B$.  It is
upper-semicontinuous from (12) and lower-semicontinuous 
from (13). 

Returning to the dilation equation, the matrices $A$ and 
$B$ share the left
eigenvector 
$(1,1,1)$.  On the complementary subspace, they reduce to 
$$A^{\prime} = \bmatrix c_{0} & 0  \\ -c_{3} & 1-c_{0} -
c_{3}\endbmatrix \ \ \ \text{and}\ \ \ B^{\prime} = \bmatrix
1-c_{0}-c_{3} & -c_{0} \\ 0 & c_{3}\endbmatrix.$$
It is $\rho(A^{\prime},B^{\prime})$ that decides the size 
of $\phi
(x) - \phi (y)$. Continuity follows from $\rho < 1$ [16]. 
Then
$\phi$ 
and $W$ belong to $C^{\alpha}$ for all $\alpha$ less than
$-\log_{2}\rho$.  (When 
$\alpha > 1$, derivatives of integer order $[\alpha]$
have H\"{o}lder exponent $\alpha - [\alpha ]$.)
In Sobolev
spaces $H^{s}$, Eirola and Villemoes
[17, 18] showed how an ordinary spectral radius ---
{\it computable} --- gives the exact regularity $s$.

\heading 5. Accuracy and orthogonality \endheading

For the Daubechies coefficients, the dilation equation 
does produce a
continuous $\phi (x)$ with H\"{o}lder exponent $0.55$ (it is
differentiable almost everywhere).  Then (8) constructs 
the wavelet.
Figure 2 shows $\phi$ and $W$ with 
$c_{0}$, $c_{1}$, $c_{2}$, $c_{3} = \frac{1}{4} (1 + 
\sqrt{3})$, 
$\frac{1}{4}(3+\sqrt{3})$, $\frac{1}{4}(3 -\sqrt{3})$,
$\frac{1}{4}(1-\sqrt{3})$.

\midspace{25.5pc}
\caption{{\smc Figure} 2.
{\rm The family $W_{4}(2^{j}x-k)$ is orthogonal.
Translates of $D_{4}$ can reproduce any $ax + b$.
Daubechies also found $D_{2p}$ with orthogonality and
$p$\<th order accuracy.}}

What is special about the four Daubechies coefficients?  
They satisfy
the requirement $\text{A}_{2}$ for second-order accuracy 
and the
separate requirement $\text{O}$ for orthogonality.  We can 
state Condition
$\text{A}_{2}$  in several forms.  In terms of
$W$, the moments $\int W(x) \, dx$ and 
$\int x \, W(x) \, dx$ are zero.  Then the Fourier 
transform of (8) yields
$P(\pi) = P^{\prime}(\pi ) = 0$.  In terms of the $c$\<'s 
(or the symbol
$P(\xi ) = \frac{1}{2} \sum c_{k} \, e^{ik\xi}$), the 
condition for
accuracy of order $p$ is $\text{A}_{p}$:  
$$\sum (-1)^{k} \, k^{m} \, c_{k} = 0  \
\text{for}\ m < p \ \ \
\text{or equivalently}\ \ \ P(\xi + \pi) = O(|\xi |^{p}). 
\tag{14}$$
This assures that translates of $\phi$ reproduce (locally) 
the powers
$1,x,\dots,x^{p-1}$ [19]. The zero moments are the 
orthogonality of
these powers to $W$.  Then the Taylor
series of $f(x)$  can be
matched to degree $p$ at each meshpoint.  The error in 
wavelet
approximation is of order $h^{p}$, where $h =
2^{-j}$ 
is the mesh width or translation step of the local functions
$W(2^{j}x)$.  The price for each extra order of accuracy 
is two extra
coefficients $c_{k}$ --- which spreads the support of 
$\phi$ and $W$
by two intervals.  A reasonable compromise is $p = 3$.
The new short wavelets may offer an alternative.

Condition $\text{A}_{p}$ also produces zeros in the 
infinite product
${\hat \phi}(\xi) = \Pi\, P(\xi /2^{j})$.  Every nonzero 
integer has the
form $n = 2^{j-1}m,\ m$ odd.  Then ${\hat \phi}(2\pi n)$ 
has the factor
$P(2 \pi n/2^{j}) = P(m\pi ) = P(\pi )$.  Therefore, the 
$p$\<th order
zero at $\xi = \pi$ in Condition $\text{A}_{p}$ ensures a 
$p$\<th order zero
of ${\hat \phi}$ at each $\xi = 2 \pi n$.  This is the 
test for the
translates of $\phi$ to reproduce $1,x,\dots,x^{p-1}$.  
That step
closes\ the circle and means approximation to order $p$.  
Please
forgive this brief recapitulation of an older theory --- 
the novelty
of wavelets is their orthogonality.  This is tested by 
Condition O:
$$\sum c_{k} \, c_{k-2m} = 2\,\delta_{0m}\ \ \ \text{or 
equivalently}\
\ \ |P(\xi)|^{2} + |P(\xi + \pi)|^{2} \equiv 1. \tag{15}$$
The first condition follows directly from $(\phi (x), \phi 
(x-m)) =
\delta_{0m}$.
The dilation equation converts this to
$(\sum c_{k} \, \phi (2x-k),\ \sum c_{\ell} \, \phi 
(2x-2m-\ell)) =
\delta_{0m}$.  It is 
the ``perfect reconstruction condition'' of digital signal
processing [20--22].  It assures that the $L^{2}$ norm is 
preserved,
when the signal $f(x)$ is separated by a low-pass filter 
$L$ and a high-pass
filter $H$. The two parts have $\parallel Lf\parallel^{2} + 
\parallel
Hf\parallel^{2} = \parallel f\parallel^{2}$.
{\it A filter is just a convolution}.
In frequency space that makes it a
multiplication.   {\it Low-pass} means that constants and 
low
frequencies survive --- we multiply by a symbol $P(\xi)$ 
that is near
$1$ for small $|\xi |$.  {\it High-pass} 
means the opposite, and for wavelets the multiplier is 
essentially
$P(\xi + \pi)$.  The two convolutions are ``mirror 
filters''.  

In the discrete case, the filters $L$ and $H$ (with 
downsampling to
remove every second row) fit into an orthogonal matrix:
$$\bmatrix \\ \ \ L \ \ \\ \\ \ \ H \ \ \\ \\ \endbmatrix =
\frac{1}{\sqrt{2}} \bmatrix
& & c_{0} & c_{1} & c_{2} & c_{3} \\
& &       &       & c_{0} & c_{1} & c_{2} & c_{3} & & \\
& & \cdot & \cdot & \cdot & \cdot & \cdot & \cdot \\
& & c_{3} & -c_{2} & c_{1} & -c_{0} \\
& &       &        & c_{3} & -c_{2} & c_{1} & -c_{0} 
\endbmatrix.
\tag 16$$
This matrix enters each step of the wavelet transform, 
from vector
$y$ to wavelet coefficients $b$.
The pyramid algorithm executes that transform by
{\it recursion with rescaling}.
We display two steps for a general wavelet and then 
specifically
for Haar on $[0,1]$:
$$\bmatrix
{\ssize L} \\ {\ssize H} \\ \vspace{1 \jot} & I & \\ \\ 
\endbmatrix
\bmatrix \\ \ \ L \ \ \\ \vspace{1 \jot} \ \ H \ \ \\ \\ 
\endbmatrix
 \text{ is } 
\frac{1}{\sqrt 2} \bmatrix
1 & 1  \\ 1 & -1  \\
& & \sqrt 2  \\ & & & \sqrt 2 \endbmatrix
\frac{1}{\sqrt 2} \bmatrix
1 & 1 \\ & & 1 & 1 \\ 1 & -1 \\ & & 1 & -1 \endbmatrix\<.
\tag 17$$
This product is still an orthogonal matrix.
When the columns of $W_4$ in \S1 are normalized to be
unit vectors, this is its inverse (and its transpose).
The recursion decomposes a function into wavelets, and the
reverse algorithm reconstructs it.
The $2$ by $2$ matrix has low-pass
coefficients $1,1$ from $\phi$ and high-pass coefficients 
$1,-1$ from
$W$.  Normalized by $\frac{1}{2}$, they satisfy Condition 
O (note
$e^{i\pi} = -1$), and they preserve the $\ell^2$ norm:
$$\left| \frac{1 + e^{i\xi}}{2}\right|^{2}\ \ + \ \ \left| 
\frac{1 +
e^{i(\xi + \pi)}}{2}\right|^{2} \equiv 1.$$
Figure~3 shows how those terms $|P(\xi)|^2$ and $|P(\xi+
\pi)|^2$
are mirror functions that add to $1$.
It also shows how four coefficients give a flatter 
response --- with
higher accuracy at $\xi = 0$.
Then $|P|^2$ has a fourth-order zero at $\xi = \pi$.

\midspace{13pc}
\caption{{\smc Figure }3.
{\rm Condition O for Haar ($p=1$) and Daubechies ($p=2$).}}

The design of filters (the choice of convolution) is a 
central
problem of signal processing --- a field of enormous size 
and
importance.  The natural formulation is in frequency 
space.  Its
application here is to multirate filters and ``subband 
coding'', with
a sequence of scales $2^{j}x$.

\rem{Note} Orthogonality of the family $\phi (x-k)$ leads by
the Poisson summation formula to $\sum
|{\hat \phi}(\xi + 2\pi n)|^{2} =1 $.  Applying the dilation
equation (7) and separating even $n$ from odd $n$ shows how
the second form of Condition O is connected to 
orthogonality:
$$\align
& \sum |{\hat \phi}(\xi + 2\pi n)|^{2} \\
&\qquad = \sum \left|P\left(\frac{\xi}{2} + \pi 
n\right)\right|^{2} \,
    \left|{\hat \phi}\left(\frac{\xi}{2} + \pi 
n\right)\right|^{2} \\
&\qquad = \left|P\left(\frac{\xi}{2}
\right)\right|^{2} \sum \left|{\hat 
\phi}\left(\frac{\xi}{2} + 
\pi 2m\right)\right|^{2}
    +  \left|P\left(\frac{\xi}{2} + \pi\right)\right|^{2}
        \sum \left|{\hat \phi}\left(\frac{\xi}{2} + \pi(2m+
1)
\right)\right|^{2} \\
&\qquad = \left|P\left(\frac{\xi}{2}\right)\right|^{2} 
+ \left|P\left(\frac{\xi}{2} + \pi\right)\right|^{2} \ \ \
    (=1\ \text{by Condition O}).
\endalign$$
The same ideas apply to $W$.  For dilation by $3^{j}$ or 
$M^{j}$ instead of
$2^{j}$, 
Heller has constructed [23] the two wavelets or $M-1$ 
wavelets that
yield approximation of order $p$.  The orthogonality 
condition 
becomes $\sum^{M-1}_{0}\ 
|P(\xi +
{2\pi j}/{M})|^{2} = 1$.

We note a technical hypothesis that must be added to 
Condition O.  It
was found by Cohen and in a new form by Lawton (see [24, 
pp.\ 
177--194]).  Without it, $c_{0} = c_{3} = 1$ passes test O.
Those coefficients give a stretched box function
$\phi = \frac{1}{3}\chi_{[0,3]}$ that is not
orthogonal to $\phi (x-1)$.  The matrix with $L$ and $H$ 
above will
be only an isometry --- it has columns of zeros.  The 
filters satisfy
$LL^{*} = HH^{*} = I$ and $LH^{*} = HL^{*} = 0$ but not 
$L^{*}L +
H^{*}H = I$.  The extra hypothesis is applied to this 
matrix $A$, or
after Fourier transform to the operator~${\Cal A}$:
$$A_{ij} = \sum^{N}_{0} c_{k} c_{j-2i+k}\ \ \ \text{or}\ \ 
\ {\Cal
A}f(\xi) = \left|P\left(\frac{\xi}{2}\right)\right|^{2}f
\left(\frac{\xi}{2}\right) \ + \
\left|P\left(\frac{\xi}{2} + \pi
\right)\right|^{2}f\left(\frac{\xi}{2} + \pi\right).$$
The matrix $A$ with $|i| < N$ and $|j| < N$ has two 
eigenvectors for
$\lambda = 1$.  Their components are $v_{m} = \delta_{0m}$ 
and $w_{m}
= (\phi (x),\phi (x-m))$.  Those must be the same!  Then 
the extra
condition, added to O, is that $\lambda = 1$ shall be a 
simple eigenvalue.
\endrem

In summary, Daubechies used the minimum number $2p$
of coefficients $c_{k}$ to satisfy the accuracy condition 
$\text{A}_{p}$
together with orthogonality.  These wavelets furnish {\it 
unconditional bases}
for the key spaces of harmonic analysis ($L^{p}$, 
H\"{o}lder, Besov,
Hardy space $H^{1}$,
$BMO$, $\dots$).  The Haar-Walsh construction fits 
functions with no extra
smoothness [25]. Higher-order wavelets fit Sobolev spaces, 
where
functions have derivatives in $L^{p}$ (see \cite {11,
pp.\ 24--27}).  With marginal exponent $p = 1$ or even $p 
< 1$, the wavelet
transform still maps onto the right discrete spaces.

\heading 6. The contest: Fourier vs. wavelets \endheading

This brief report is included to give some idea of the 
decisions now
being reached about standards for video compression.  The 
reader will
understand that the practical and financial consequences 
are very
great.  Starting from an image in which each color at each 
small square
(pixel) is
assigned a numerical shading between $0$ and $255$, the 
goal is to
compress all that data to reduce the transmission cost.  
Since $256 = 2^{8}$,
we have $8$ bits for each of red-green-blue.  The
bit-rate of transmission is set by the channel capacity, the
compression rule is decided by the filters and quantizers, 
and the picture
quality 
is subjective.  Standard images are so familiar that 
experts know
what to look for --- like tasting wine or tea.

Think of the problem mathematically.  We are given 
$f(x,y,t)$, with
$x$-$y$ axes on the TV screen and the image $f$ changing 
with time $t$.
 For digital signals all variables are discrete, but a 
continuous
function is close --- or piecewise continuous when the 
image has
edges.  Probably $f$ changes gradually as the camera 
moves.  We could
treat $f$ as a sequence of still images to compress 
independently,
which seems inefficient.  But the direction of movement is
unpredictable, and too much effort spent on extrapolation 
is also inefficient.
A compromise is to encode every fifth or tenth image and, 
between
those, to work with the time differences $\Delta f$ --- 
which have less
information and can be compressed further.

Fourier methods generally use real transforms (cosines). 
The picture is broken
into blocks, often 8 by 8. This improvement in the scale 
length is more
important than the control of $\log n$ in the FFT cost. 
(It may well be more
important than the choice of Fourier.)
After twenty years of refinement,
the algorithms are still being fought over and  improved.  
Wavelets
are a recent entry, not yet among the
heavyweights.  The accuracy test $\text{A}_{p}$ is often 
set aside in the
goal of constructing ``brick wall filters'' --- whose 
symbols
$P(\xi)$ are near to characteristic functions.  An exact 
zero-one
function in Figure~3 is of course impossible --- the 
designers are
frustrated by a small theorem in mathematics.
(Compact support of $f$ and ${\hat f}$ occurs only for $f 
\equiv 0$.)
In any case the Fourier transform of a step function 
has oscillations that can murder 
a pleasing signal --- so a compromise is reached.

Orthogonality is not set aside.  It is the key constraint. 
 There may
be eight or more bands ($8$ times $8$ in two dimensions) 
instead of two.
Condition O has at least eight terms $|P(\xi + k\pi 
/8)|^{2}$.
After applying the convolutions, the energy or entropy in
the high frequencies is usually small
and the compression of that part of the signal is 
increased --- to
avoid wasting bits.  The actual
encoding or ``quantization'' is a separate and very subtle 
problem,
mapping the 
real numbers to $\{ 1,\dots ,N\}$.  A vector quantizer is 
a map from $\bold
R^d$,
and the best are not just tensor products [28].  Its 
construction is
probably more important to a successful compression than 
refining the
filter.

Audio signals have fewer dimensions and more bands --- as 
many as $512$.  One
goal of compression is a smaller CD disk.  Auditory 
information seems
to come in 
octaves of roughly equal energy --- the energy density 
decays
like $1/\xi$.  Also physically, the cochlea has several 
critical bands per
octave.  (An active problem in audio compression is to use 
psychoacoustic
information about the ear.)  Since $\int d\xi /\xi$ is the 
same from $1$ to $2$ and
$2$ to $4$ and $4$ to $8$ (by a theorem we teach 
freshmen!), subband
coding stands a good chance. 

That is a barely adequate description of a fascinating 
contest. It is
applied analysis (and maybe free enterprise) at its best.  
For video
compression, the Motion Picture Experts Group held a 
competition in
Japan late in 1991.  About thirty companies entered 
algorithms.  Most were
based on cosine transforms, a few on wavelets.  The best
were all windowed Fourier.  Wavelets
were down the list but not unhappy.
Consolation was freely offered and accepted.  The choice 
for HDTV,
with high definition, may be 
different from this MPEG standard to send a rougher 
picture at a
lower bit-rate.

{\it I must emphasize}\/:  
The real contest is far from over.  There are promising
wavelets (Wilson bases and coiflets) that were too recent 
to enter.  Hardware
is only beginning to come---the first wavelet chips are 
available.  MPEG did
not see the best that all transforms can do.

In principle, wavelets are better for images, and Fourier 
is the right choice
for music.  Images have sharp edges; music is sinusoidal.  
The $j$\<th Fourier 
coefficient of a step function is of order $1/j$.  The 
wavelet coefficients
(mostly zero) are multiples of $2^{-j/2}$.  The $L^2$ 
error drops
exponentially,
not polynomially, when $N$ terms are kept.  To confirm 
this comparison, Donoho
took digitized photos of his statistics class. He 
discarded $95\%$ of the
wavelet and the Fourier coefficients, kept the largest 
$5\%$, and reconstructed
two pictures.  (The wavelets were ``coiflets'' \cite {24}, 
with greater
smoothness and symmetry but longer support.  Fourier 
blocks were not tried.)
Every student preferred the picture from wavelets.

The underlying rule for basis functions seems to be this:  
choose scale lengths
that match the image and allow for spatial variability.  
Smoothness is visually
important, and $D_4$ is being superseded.  Wavelets are 
not the only possible
construction, but they have opened  the door to new bases. 
 In the mathematical
contest (perhaps eventually in the business contest) 
unconditional bases are the
winners. 

We close by mentioning fingerprints.  The FBI has more 
than 30
million in filing cabinets, counting only criminals.  
Comparing one
to thousands of others is a daunting task.  Every 
improvement leads
to new matches and the solution of old crimes.  The images 
need to
be digitized.

The definitive information for matching fingerprints is in 
the ``minutiae'' of
ridge endings and bifurcations [29].  At 500 pixels per 
inch, with
256 levels of gray, each card has $10^{7}$ bytes of data.
Compression is essential and $20:1$ is the goal.  The 
standard from
the Joint Photographic Experts Group (JPEG) is 
Fourier-based,
with $8$ by $8$ blocks, and the ridges are broken.  The 
competition
is now between wavelet algorithms associated with Los 
Alamos and Yale
[30--33] --- fixed basis versus ``best basis'', $\ell < 
100$ subbands
or $\ell > 1000$, vector or scalar quantization.  There is 
also a
choice of coding for wavelet coefficients (mostly near 
zero when the
basis is good). The best wavelets may be {\it 
biorthogonal}  ---
coming from two wavelets $W_{1}$ and $W_{2}$.  This allows a
left-right symmetry [24], which is absent in Figure~2.  The
fingerprint decision is a true contest in applying pure 
mathematics.

\heading Acknowledgment \endheading

I thank Peter Heller for a long conversation about the 
MPEG contest and its
rules.

\rem{Additional note} After completing this paper I
learned, with pleasure and amazement, that a thesis which 
I had promised 
to supervise (``formally'', in the most informal sense
of that word) was to contain the filter design for MIT's 
entry in the
HDTV competition.  The Ph.D. candidate is Peter Monta.  The
competition is still ahead (in 1992).  Whether won or 
lost, I am sure
the degree will be granted!  These paragraphs briefly 
indicate how
the standards for High Definition Television aim to yield 
a very
sharp picture.

The key is high resolution, which requires a higher 
bit-rate of
transmission.  For the MPEG contest in Japan --- to 
compress videos onto CD's
and computers --- the rate was $1$ megabit/second.  For 
the HDTV contest
that number is closer to $24$.  Both compression ratios 
are about
$100$ to $1$.  (The better picture has more pixels.)  The 
audio
signal gets $\frac{1}{2}$ megabit/sec for its four stereo 
channels;
closed captions use less.  In contrast,
conventional television has no compression at all --- in 
principle, 
you see everything.  The color standard was set in 1953, 
and the black
and white standard about 1941.

The FCC will judge between an AT\&T/Zenith entry, two 
MIT/General
Instruments entries, and a partly European entry from 
Philips and
others.  These 
finalists are all digital, an advance which surprised the 
{\it New York
Times\/}.  Monta proposed a filter that uses seven 
coefficients or
``taps'' for low-pass and four for high-pass.  Thus the 
filters are
not mirror images as in wavelets, or brick walls either.
Two-dimensional images come from tensor products of 
one-dimensional
filters.  Their exact coefficients will not be set until 
the last
minute, possibly for secrecy --- and cosine transforms may 
still be
chosen in the end. 

The red-green-blue components are converted by a $3$ by $3$
orthogonal matrix to better coordinates.  Linear algebra 
enters,
literally the spectral theorem.  The luminance axis from the
leading eigenvector gives the brightness.

A critical step is motion estimation, to give a quick and 
close
prediction of successive images.  A motion vector is 
estimated for
each region in the image [34].  The system transmits only 
the
difference between predicted and actual images --- the 
``motion
compensated residual''.  When that has too much energy, 
the motion
estimator is disabled and the most recent image is sent.  
This will be
the case when there is a scene change.  Note that coding 
decisions
are based on the {\it energy} in
different bands (the size of Fourier coefficients).  The 
$L^{1}$ norm is
probably better.  Other features may be
used in 2001.

It is very impressive to see an HDTV image.  The final 
verdict has just
been promised for the spring of 1993.  Wavelets will not 
be in that
standard, but they have no shortage of potential 
applications [24,
35--37]. A recent one is the LANDSAT 8 satellite, which 
will locate a grid on
the earth with pixel width of 2 yards.  The compression 
algorithm that does
win will use good mathematics.
\endrem

\enddocument